\documentclass[12pt,draft]{amsart}
\usepackage{amsmath,amsthm,latexsym,amscd,amsbsy,amssymb}
 \relax

\chardef\bslash=`\\ 

\makeatletter
\def\verbatim{\interlinepenalty\@M \@verbatim
  \leftskip\@totalleftmargin\advance\leftskip2pc
  \frenchspacing\@vobeyspaces \@xverbatim}
\makeatother
\makeatletter
  \def\dgt@k{\dg@DX=1 \dg@DY=-4 \dg@SIZE=1}
\makeatother

\theoremstyle{plain}
\newtheorem{thm}{Theorem}[section]
\newtheorem{cor}[thm]{Corollary}

\newtheorem{pro}[thm]{Proposition}

\theoremstyle{definition}

\newcounter{rmnum}

\def\ep{\varepsilon}

\numberwithin{equation}{section}

\begin{document}


\title[Real rank and squaring mapping for unital $C^{\ast}$-algebras]
{Real rank and squaring mapping for unital $C^{\ast}$-algebras}
\author{A.~Chigogidze}
\address{Department of Mathematics and Statistics,
University of Saskatche\-wan,
McLean Hall, 106 Wiggins Road, Saskatoon, SK, S7N 5E6,
Canada}
\email{chigogid@math.usask.ca}
\thanks{The first named author was partially supported by NSERC research grant.}
\author{A. Karasev}
\address{Department of Mathematics and Statistics,
University of Saskatche\-wan,
McLean Hall, 106 Wiggins Road, Saskatoon, SK, S7N 5E6,
Canada}
\email{karasev@math.usask.ca}
\author{M.~R{\O}rdam}
\address{Department of Mathematics,
University of Copenhagen,
Universitetsparken~5,
2100 Copenhagen {\O},
Denmark}
\email{rordam@math.ku.dk}



\begin{abstract}{It is proved that if $X$ is a compact Hausdorff space
    of Lebesgue dimension $\dim(X)$, then the squaring
    mapping $\alpha_{m} \colon \left( 
      C(X)_{\mathrm{sa}}\right)^{m} \to C(X)_{+}$, defined by
    $\alpha_{m}(f_{1},\dots ,f_{m}) = \sum_{i=1}^{m} f_{i}^{2}$, is
    open if and only if $m -1 \ge \dim(X)$. Hence the Lebesgue
    dimension of $X$ can be detected from openness of the squaring
    maps $\alpha_m$. In the case $m=1$ it is proved that the map $x
    \mapsto x^2$, from the self-adjoint elements of a unital
    $C^{\ast}$-algebra $A$ into its positive elements, is open if and
    only if $A$ is isomorphic to $C(X)$ for some compact Hausdorff
    space $X$ with $\dim(X)=0$.
} 
\end{abstract}

\maketitle \markboth{A.~Chigogidze, A.~Karasev, M.~R{\o}rdam}{Real
  rank and squaring mapping for unital $C^{\ast}$-algebras} 

\section{Introduction}\label{S:intro}
A compact Hausdorff space $X$ is defined to have Lebesgue dimension
$\le m$ if for every closed subset $F$ of $X$, each continuous map $F
\to S^m$ has a continuous extension $X \to S^m$. 

Various types of ranks for (unital) $C^{\ast}$-algebras
have been inspired by corresponding prototypes in the classical
dimension theory of (compact) spaces, such as the one given
above. While the Lebesgue dimension of a 
compact space has numerous equivalent formulations, the extensions of
these equivalent formulations to non-commutative $C^{\ast}$-algebras
most often differ. 
Examples of such ranks for
$C^{\ast}$-algebras are 
the \emph{stable rank} defined by Rieffel in \cite{rieffel}, the  
\emph{real rank} defined by Brown and Pedersen in \cite{brownped91},
the \emph{analytic rank} defined by Murphy in \cite{murphy2}, the
\emph{completely positive rank} considered by Winter in \cite{winter},
and the \emph{bounded rank} defined in \cite{chiva}.
Other ranks of 
$C^{\ast}$-algebras are the \emph{tracial rank} defined by Lin in
\cite{lin} and the \emph{exponential rank} defined by Phillips in
\cite{phillips}. 

It was shown in \cite{chiva} that a 
unital $C^{\ast}$-algebra $A$ has real rank at most $n$ if the
squaring map $(x_{1},\dots ,x_{n+1}) \mapsto
\sum_{i=1}^{n+1}x_{i}^{2}$, from the set of $(n+1)$-tuples of
self-adjoint elements to the set of positive elements in $A$, is
open; and it was asked if the reverse also holds, in which case
openness of the squaring maps would determine the real rank of the
$C^{\ast}$-algebra.

In the present note we answer this question 
in the affirmative in the commutative case --- and in the negative in
the general (non-commutative) case. The latter --- negative --- answer
follows from our result that the squaring map $x \mapsto x^2$ (from
the set of self-adjoint elements to the set of positive elements) is
open if and only if the $C^{\ast}$-algebra is commutative and of real
rank zero. Hence the squaring map $x \mapsto x^2$ is not open for the
$C^{\ast}$-algebra, $M_2$, of $2$ by $2$ matrices, but this
$C^{\ast}$-algebra has real rank zero.



\section{Results}\label{S:results}
For a $C^{\ast}$-algebra $A$ we use the standard notation
$A_{\mathrm{sa}}$ and $A_{+}$ to denote the set of all self-adjoint
and the set of all positive elements of $A$, respectively. 
The real rank of a unital $C^{\ast}$-algebra $A$, denoted by
$\operatorname{RR}(A)$, 
is in \cite{brownped91} defined as follows: 
For each non-negative integer $n$, $\operatorname{RR}(A) \leq n$ if for every
$(n+1)$-tuple $(x_{1},\dots ,x_{n+1})$ 
in $A_{\mathrm{sa}}$ and every $\varepsilon > 0$,
there exists an $(n+1)$-tuple $(y_{1},\dots ,y_{n+1})$ in $A_{\mathrm{sa}}$
such that $\sum_{k=1}^{n+1} y_{k}^{2}$ is invertible and
$\sum_{k=1}^{n+1}\|x_k-y_k\| < \varepsilon$.  

Let us say that a unital $C^{\ast}$-algebra $A$ has  an open
$m$-squaring map if the map $\alpha_{m} \colon \left(
 A_{sa}\right)^{m} \to A_{+}$, defined by 
$\alpha_{m}(x_{1},\dots ,x_{m}) =
\sum_{k=1}^{m}x_{k}^{2}$, is open. 
Observe that $\alpha_m$ is open at $(x_1, \dots, x_m)$ if for every
$\varepsilon >0$ there is $\delta>0$ such that for all $a \in A_+$ with
$\|\sum_{k=1}^m x_k^2 - a\| < \delta$ there is an $m$-tuple $(y_1,
\dots, y_m)$ in $A_{\mathrm{sa}}$ with $\sum_{k=1}^m y_k^2 = a$ and
$\sum_{k=1}^m \|x_k-y_k\| < \ep$.

For the readers convenience we
present a shorter proof of \cite[Proposition 7.1]{chiva}. 

\begin{pro}\label{P:square}
Let $A$ be a unital $C^*$-algebra. If the $(n+1)$-squaring map on $A$
is open, then $\operatorname{RR}(A) \leq n$.
\end{pro}
\begin{proof} Let $(x_1, \dots, x_{n+1})$ be an $(n+1)$-tuple of
  self-adjoint elements in $A$ and let $\varepsilon >0$. By openness
  of the $(n+1)$-squaring map there is $\delta >0$ and an
  $(n+1)$-tuple $(y_1,\dots, y_{n+1})$ of self-adjoint elements in $A$
  such that $\sum_{k=1}^{n+1}\| x_k - y_k \| < \varepsilon$ and 
  $\sum_{k=1}^{n+1} y_k^2 = \sum_{k=1}^{n+1} x_k^2 + \delta  \cdot
   1$, and the latter element is invertible (because each $x_k^2$ is
  positive). 
\end{proof}

\noindent Next, we prove the reverse of Proposition~\ref{P:square} in the
commutative case.


\begin{pro}\label{P:abelian}
If $X$ is a compact space such that $\dim X \leq n$, then $C(X)$ has open $(n+1)$-squaring map.
\end{pro}
\begin{proof}
Let $f = f_{1}^{2} +\dots +f_{n+1}^{2}$. 
Put $m_i =\sup f_i$, $i=1,\dots ,n+1$ and let $m=\max\{ m_i\}$.
Fix $\ep >0$ and let 
\[ \delta =\min \left\{\frac{ (\ep
    /3)^4}{m^2},\left(\frac{\ep}{3}\right)^2 \right\}\;\text{and}\; 
U=\left\{ x\in X\colon f(x)> \left(\frac{\ep}{3}\right)^2 \right\} .\]

Let also 
$A = f^{-1}([0, (\varepsilon/3)^2])$ and $S =
f^{-1}((\varepsilon/3)^2)$. 
Then $A$ and $S$ are closed subsets of $X$ such that $A= X\backslash
U$ and $S\subseteq A$. 

Now consider the diagonal product

\[ F(x) = \triangle\{ f_{1}(x), \dots , f_{n+1}(x)\} \colon X \to {\mathbb R}^{n+1}\]

\noindent and note that

\[ A = F^{-1}\left( B^{n+1}\right) \;\text{and}\; S =F^{-1}\left(
  S^{n}\right) , \] 

\noindent where 
\[ B^{n+1} =
B^{n+1}\left(\mathbf{0},\ep/3\right)\;\text{and}\; S^{n} =
\partial B^{n+1}\left(\mathbf{0},\ep/3\right) \]

Since $\dim A\leq \dim X \leq n$, the map $F|S \colon S \to S^{n}$
admits an extension $H \colon A \to S^{n}$ (see, for instance,
\cite[Ch.\ 3, Theorem 2.2]{pears}).  Represent $H$ as the diagonal
product $H = \triangle\{ h_{i}, \dots ,h_{n+1}\}$ , where each $h_{i}$
maps $A$ into ${\mathbb R}$. Since $H(A) \subseteq S^{n}$ it follows
that $ h_{1}^{2} +\dots +h_{n+1}^{2} = (\ep/3)^2$. Note also that
since $H|S = F|S$ we 
have $h_{i}|S = f_{i}|S$ for each $i=1,\dots ,n+1$. 

The last condition allows us to define for each $i= 1,\dots ,n+1$ a
continuous map $\widetilde{h_i}$ on $X$ by letting 
\[
 \widetilde{h_i}(x)=
\begin{cases}
 f_i (x), &\text{if} \; x\in U \\
       h_{i}(x), &\text{if} \; x\in A
\end{cases}
\]

Observe that the function $\widetilde{h} = \widetilde{h_1}^2 +\dots
+\widetilde{h}_{n+1}^2$ is strictly positive on $X$. Notice also that 
$\widetilde{h} |_U =f|_U$ and $\widetilde{h} |_A =(\ep/3)^2$.

Next consider $g\in O(f,\delta )\cap C_{+}(X)$ and define a function
$\lambda$ on $X$ by 
$\lambda (x)= \big(g(x)/\widetilde{h}(x)\big)^{1/2}$. 
Note that $\lambda\geq 0$.

Now define $g_i$, $i=1,\dots ,n+1$, on $X$ by the formula $g_{i}(x)
=\widetilde{h}_{i}(x)\cdot \lambda(x)$.  
Clearly 

\begin{multline*}
 g_{1}^{2}(x) +\dots +g_{n+1}^{2}(x)  = \left(\widetilde{h}_{1}(x)\cdot \lambda(x)\right)^{2} +\dots +\left(\widetilde{h}_{n+1}(x)\cdot \lambda(x)\right)^{2} =\\
 \lambda^{2}(x)\left( \widetilde{h}_{1}^{2}(x)+\dots +\widetilde{h}_{n+1}^{2}(x)\right) = \frac{g(x)}{\widetilde{h}(x)}\cdot\left( \widetilde{h}_{1}^{2}(x)+\dots +\widetilde{h}_{n+1}^{2}(x)\right) = g(x) .
\end{multline*}

Next let us show that $g_{i}$ is sufficiently close to $f_{i}$ for each $i = 1,\dots ,n+1$. Indeed, since $g\in O(f,\delta )$ we have for each $x \in A$ and therefore 
\[ g(x) < f(x) +\delta < \left(\frac{\ep}{3}\right)^{2}+\left(\frac{\ep}{3}\right)^{2}.\]

\noindent Since $g_i^2\le g$ and $f_i^2\le f$, the last inequality implies that

\[ |g_i (x)|<\sqrt{2\left(\frac{\ep}{3}\right) ^2}< 2\frac{\ep}{3}\mbox{ and }|f_i (x)|\le\frac{\ep}{3} \]

\noindent for all $x\in A$.
Hence 
\[ |f_i (x)-g_i (x)|\le |f_i (x)|+|g_i (x)|<\frac{\ep}{3} +2\frac{\ep}{3} =\ep \]
as required.

Further, if $x \in U$, then $\widetilde{h}(x)=f(x)$ and consequently 
\[ \widetilde{h}(x)-\delta <g(x)<\widetilde{h}(x)+\delta .\]
Hence
\[ 1-\frac{\delta}{\widetilde{h}(x)}<\frac{g(x)}{\widetilde{h}(x)}<1+\frac{\delta}{\widetilde{h}(x)} \]

\noindent for $x \in U$.
Since
\[ \widetilde{h}(x)=f(x) > \left(\frac{\ep}{3} \right) ^2\text{for}\; x \in U\;\text{and}\;\delta\le\frac{(\ep /3)^4}{m^2} \]
we have (for $x \in U$)

\[
 1-\left(\frac{\ep}{3m}\right) ^2 = 1 -\frac{\frac{1}{m^2}\cdot\left( \frac{\ep}{3}\right)^{4}}{\left(\frac{\ep}{3}\right)^2} < \frac{g(x)}{\widetilde{h}(x)}
< 1 +\frac{\frac{1}{m^2}\cdot\left( \frac{\ep}{3}\right)^{4}}{\left(\frac{\ep}{3}\right)^2} = 1+\left(\frac{\ep}{3m}\right) ^2  \]
\noindent and

\[ 1-\left(\frac{\ep}{3m}\right) ^2 < \lambda^{2}(x) < 1+\left(\frac{\ep}{3m}\right) ^2 .\] 
Consequently,
\[ \left| 1-\lambda^{2}(x)\right| < \left(\frac{\ep}{3m}\right)^{2}.\]

\noindent Since $\lambda (x) \geq 0$, this implies
\[\left[1-\lambda (x)\right]^{2} \leq |1-\lambda (x)|\cdot |1+\lambda (x)| = | 1-\lambda^{2}(x)| < \left(\frac{\ep}{3m}\right)^{2} .\]

\noindent Therefore 

\[ \left| 1-\lambda (x)\right| \leq \frac{\ep}{3m}\;\;\text{for any}\;\; x \in U .\]

Finally we have
\[ |f_i(x)-g_i(x)|=|1-\lambda (x)|\cdot |f_i(x)| <\frac{\ep}{3m}\cdot m <\frac{\ep}{3}\;\text{for any}\; x\in U .\]

This completes the verification of the fact that $\left| f_{i}(x) - g_{i}(x)\right| < \ep$ for each $x \in X$ and any $i=1,\dots ,n+1$.
\end{proof}

\begin{cor}\label{C:zero}
Let $A$ be a unital $C^{\ast}$-algebra. Then the following conditions are
equivalent:
\begin{itemize}
\item[(i)]
The squaring map $x \mapsto x^2$ from $A_\mathrm{sa}$ to $A_+$ is
open. 
\item[(ii)] $A$ is commutative and $\operatorname{RR}(A) = 0$.
\item[(iii)]
$A$ is  isomorphic to a $C^{\ast}$-algebra of the form $C(X)$ for a
compact Hausdorff space $X$ with $\dim X =0$. 
\end{itemize}
\end{cor}
\begin{proof} The equivalence of (ii) and (iii) follows from Gelfand's
  duality and \cite[Proposition 1.1]{brownped91}.

The implication (iii) $\Rightarrow$ (i) follows from Proposition
\ref{P:abelian}. 

(i) $\Rightarrow$ (ii). Assume that (i)
holds. Then $\operatorname{RR}(A) = 0$ by Proposition
\ref{P:square}. It remains to show that $A$ is 
commutative. Since $A$ is of real rank zero it suffices to show that
any two projections $p,q$ in $A$ commute.

Take the symmetry $s = p-(1-p)$. Then $s$ is self-adjoint and
$s^2=1$. By openness of the squaring map there are self-adjoint
elements $s_n$ in $A$ such that $\|s_n-s\| \to 0$ and $s_n^2 =
1+n^{-1}q$. Define $\varphi \colon {\mathbb R} \to {\mathbb R}$ by
$\varphi(t) = \max\{0,t\}$. For each $n$, the element $\varphi(s_n)$
commutes with $s_n$, hence with $s_n^2$, and hence with $q$. Since
$\varphi(s)=p$ we obtain
$$pq-qp = \lim_{n \to \infty} (\varphi(s_n)q - q\varphi(s_n)) = 0,$$
as desired.
\end{proof}

\section{Related comments and open problems}

\noindent 
\emph{Existence of square roots:} Suppose that $A$ is a unital
$C^{\ast}$-algebra and that $x$ is a self-adjoint element in $A$. Does
there exist a continuous square root $\rho_x = \rho \colon \Omega \to
A_{{\mathrm{sa}}}$ 
(i.e., $\rho(a)^2=a$ for all $a \in \Omega$)
defined on an open neighborhood $\Omega \subseteq A_+$ of $x^2$ such that
$\rho(x^2)=x$? If this is true for all self-adjoint elements $a$ in
$A$, then the equivalent conditions of Corollary~\ref{C:zero} are
satisfied. 

Suppose that $A = C(X)$ for some 0-dimensional compact Hausdorff space
$X$ (i.e., that the conditions of Corollary~\ref{C:zero} are
satisfied). Take a self-adjoint (i.e., real valued) $f
\in C(X)$, and suppose that there is a
clopen set $U$ such that $f(x) \ge 0$ for all $x \in U$ and $f(x) \le
0$ for all $x \in X \backslash U$. Then the function $\rho_U \colon
C(X)_+ \to C(X)_{{\mathrm{sa}}}$ defined by
$$\rho_U(g) = \begin{cases} \sqrt{g(x)}, & x \in U, \\ -\sqrt{g(x)}, & x
  \in X \backslash U, \end{cases}$$
is a continuous square root with $\rho_U(f^2) = f$. It is not clear to
the authors if there are continuous square roots at arbitrary real
valued functions $f$ in $C(X)$. 

In the case where $A = M_n$, the $C^{\ast}$-algebra of $n$ by $n$
matrices, if $x$ is a self-adjoint element and if $x^2$ has $n$
distinct eigenvalues, then there is a continuous square root $\rho$
with $\rho(x^2)=x$ defined on some neighborhood of $x^2$.  

In the case where $A=M_2$, it follows from Corollary~\ref{C:zero} (and its
proof) that there is 
no continuous square root $\rho$ defined on a neighborhood of $I$ such
that $\rho(I) = {\mathrm{diag}}(1,-1)$. It is easily checked
explicitly that if $r$ is a (small) non-zero real number, then any
square root of 
${{\left(\begin{smallmatrix} 1 & r \\ r & 1 \end{smallmatrix} \right)}}$
is of the form
${{\left(\begin{smallmatrix} a & s \\ s & a \end{smallmatrix} \right)}}$,
where $a$ and $s$ are real numbers satisfying $a^2 + s^2 = 1$ and
$2as=r$, and any such square root has 
distance at least 1 to ${\mathrm{diag}}(1,-1)$.

\vspace{.3cm} \noindent We end this note by listing some open problems
related to openness of the squaring maps:

\vspace{.3cm} \noindent{\bf Question 1.} Let $A$ be a unital
$C^{\ast}$-algebra, let $m$ be a positive integer, and suppose that
the squaring map $\alpha_m$ (defined above Proposition~\ref{P:square}) is
open. Does it follow that $\alpha_n$ is open for all $n \ge m$?

\vspace{.3cm} \noindent The answer to Question~1 is affirmative for
commutative $C^{\ast}$-algebras by Propositions~\ref{P:square} and
\ref{P:abelian}. The difficulty in this question lies in the fact that
if $\Omega$ is an open subset of $A_+$ and if $a \in A_+$, then
$a+\Omega$ need not be open in
$A_+$. (For instance, $1+A_+$ is not open in
$A_+$.)

\vspace{.3cm} \noindent{\bf Question 2.} Are Propositions~\ref{P:square} and
\ref{P:abelian} valid also in the \emph{non-unital} case? (For
Proposition~\ref{P:abelian}, this means that we will be talking about
locally compact Hausdorff spaces rather than compact Hausdorff
spaces.) What is the relationship between openness of $\alpha_n$ on a
non-unital $C^{\ast}$-algebra $A$ and openness of $\alpha_n$ on its
unitization?

\vspace{.3cm} \noindent{\bf Question 3.} Are the squaring maps
$\alpha_m$ open for all $m \ge 2$ when $A$ is a unital $C^{\ast}$-algebra of
real rank zero? 

\vspace{.3cm} \noindent{\bf Question 4.} Does the class 
of $C^{\ast}$-algebras, for which the squaring map $\alpha_2$ is open,
have any nice properties?
More generally, are there any justifications for considering the rank
of a $C^{\ast}$-algebra defined by openness of the squaring maps; and
will this rank reflect any ``dimension like'' properties of the
$C^{\ast}$-algebra?


\begin{thebibliography}{99}
\bibitem{brownped91}
L.~G.~Brown, G.~K.~Pedersen, {\em $C^{\ast}$-algebra
of real rank zero},
J. Functional Anal. {\bf 99} (1991), 131-149.


\bibitem{chiva}
A.~Chigogidze, V.~Valov, {\em Bounded rank of $C^{\ast}$-algebras},
preprint, 2001. 

\bibitem{lin}
H.~Lin, {\em The tracial topological rank of $C^{\ast}$-algebras},
Proc. London Math. Soc. {\bf 83} (2001), 199--234.


\bibitem{murphy2}
G.~J.~Murphy, {\em The analytic rank of a $C^{\ast}$-algebra}, Proc. Amer. Math. Soc. {\bf 115} (1992), 741--746.

\bibitem{pears}
A.~R.~Pears, {\em Dimension Theory of General Spaces}, Cambridge University Press, Cambridge, 1975.

\bibitem{phillips}
N.~C.~Phillips, {\em Simple $C^{\ast}$-algebras with the property (FU)}, Math. Scand. {\bf 69} (1991), 121--151.

\bibitem{rieffel}
M.~A.~Rieffel, {\em Dimension and stable rank in the $K$-theory of $C^{\ast}$-algebras},
Proc. London Math. Soc. {\bf 46} (1983), 301--333.

\bibitem{winter}
W.~Winter, {\em Covering dimension for nuclear $C^{\ast}$-algebras}, Preprint math.OA/0107218 (2001).

\end{thebibliography}
\end{document}